\documentclass[a4paper,10pt]{article}

\usepackage{amssymb}

\newtheorem{proposition}{Proposition}
\newtheorem{corollary}{Corollary}

\begin{document}

%\begin{frontmatter}

\title{On the  multidimensional permanent  and $q$-ary
designs\footnote{The work is supported by RFBR (grant 13-01-00463).
}}

\author{Vladimir N. Potapov
}

\maketitle

{\footnotesize Sobolev Institute of Mathematics, 4 Acad. Koptug
ave., Novosibirsk, Russia, 630090

 Novosibirsk State University, 2 Pirogova st., Novosibirsk, Russia,
630090

             email: vpotapov@math.nsc.ru}

\begin{abstract}
An $H(n,q,w,t)$ design is  a collection of some $(n-w)$-faces of the
hypercube $Q^n_q$ that perfectly pierce all $(n-t)$-faces $(n\geq
w>t)$. An $A(n,q,w,t)$ design is a collection of some $(n-t)$-faces
of $Q^n_q$ that perfectly cover all $(n-w)$-faces. The numbers of H-
and A-designs are expressed in terms of the multidimensional
permanent. Several constructions of H- and A-designs are given and
the existence of $H(2^{t+1},s2^t,2^{t+1}-1,2^{t+1}-2)$ designs is
proven for all $s,t\geq 1$.

Keywords: perfect matching; clique matching; permanent; MDS code;
Steiner system; H-design

MSC2010 05B05; 05C65

\end{abstract}

\section{Introduction}
\label{intro}  The H-design  by Hanani \cite{41:Han} is a
generalization of a Steiner system or t-design. The notation of
H-design is due to Mills \cite{41:Mil}. Let $X$ be a set of points
and let $C=\{C_1,\dots,C_n\}$ be a partition of $X$ into $n$ sets of
cardinality $q$. A {\it transverse} of $C$ is a subset of $X$
meeting each set $C_i$ at most in  one point. The  set of
$w$-element transverses of $C$ is an $H(n,q,w,t)$ {\it design}
(briefly, H-design) if each $t$-element transverse of $C$ lies  in
exactly one transverse of the H-design. We propose another
generalization of  t-design. A set of $t$-element transverses of $C$
is an $A(n,q,w,t)$ {\it design} (briefly, A-design) if each
$w$-element transverse of $C$ contains exactly one transverse of the
A-design. We imply  everywhere that $n\geq w> t\geq 1$, $q\geq 1$
and all these numbers are integer. The idea of considering the
A-designs belongs to S. V. Avgustinovich.

Put $Q_q=\{0,1,\dots,q-1\}$ and $Q_{q*}= Q_q\cup\{*\}$. It is clear
that each $w$-transverse of $C$ corresponds to the {\it codeword}
$(a_1,\dots,*,\dots,a_i,\dots,*,\dots,a_n)\in Q_{q*}^n$ where $a_i$
is the label of the element of $C_i$ that belongs to the
$w$-transverse. Position $j$ of the codeword contains $*$ if and
only if the $w$-transverse does not intersect  $C_j$. Define the
{\it weight} of a codeword from $Q_{q*}^n$ as $n$ minus the number
of symbols $*$ contained in the codeword. Then the set $H$
consisting of some vectors $x\in Q_{q*}^n$ of weight $w$  is an
$H(n,q,w,t)$ design if each vector $y\in Q_{q*}^n$ of weight $t$ is
covered by exactly one $x\in H$. Analogously the set $A$ consisting
of vectors $y\in Q_{q*}^n$ of weight $t$ is an $A(n,q,w,t)$ design
if each $x\in Q_{q*}^n$ of weight $w$ covers  exactly one  $y\in A$.

If $q=1$ then an $H(n,1,w,t)$ design is just a Steiner system
$S(t,w,n)$ (here $*$ is replaced by $0$ and  $0$ is replaced by $1$)
and an $A(n,1,w,t)$ design is just a Steiner system $S(n-w,n-t,n)$
(here  $*$ is replaced by $1$). In \cite{41:ZR} an H-design was
called a $q$-ary Steiner system. A set $T$ of  $y\in Q_{1*}^n$ of
weight $t$ is an $(n,w,t)$-Turan system, if each  $x\in Q_{1*}^n$ of
weight $w$ covers at least one  $y\in T$. Hence an $A(n,1,w,t)$
design is a special case of a $(n,w,t)$-Turan system.

The set $Q_q^n$ is called the {\it hypercube}. The set of faces of
$Q^n_q$ is in one-to-one correspondence with $Q_{q*}^n$ and each
$k$-dimensional face ($k$-face) corresponds to the codeword with $k$
symbols $*$. Thus an $H(n,q,w,t)$ design is  a~piercing consisting
of $(n-w)$-faces of  $Q^n_q$ with the property that each
$(n-t)$-face contains exactly one $(n-w)$-face of an H-design; and
an $A(n,q,w,t)$ design is a~covering consisting of $(n-t)$-faces of
$Q^n_q$ with the property that each $(n-w)$-face is contained in
exactly one $(n-t)$-face of an A-design.

If $w=n$ then an $H(n,q,w,t)$ design is just an MDS code in $Q^n_q$
with code distance $d=n-t+1$. If $w=n$ and $t=n-1$ then an
$A(n,q,w,t)$ design is just a tiling of the hypercube by 1-faces. If
$q=2$ then this tiling is equivalent to a perfect
matching\footnote{\,Here we consider $Q_2^n$ as a minimal Hamming
distance graph.} in $Q_2^n$. If $q>2$ then $A(n,q,n,n-1)$ design is
called a~{\it perfect clique matching} (see \cite{41:Pot}) because
the 1-faces of $Q_q^n$ one-to-one correspond to the maximal cliques
in the hypercube.  It is clear that $H(n,q,n,n-1)$ and $A(n,q,n,t)$
designs exist for all $q\geq 2$ and $n\geq 2$. A set of 1-faces is
called a {\it precise clique matching} if it is both
$H(n,q,n-1,n-2)$ design and $A(n,q,n,n-1)$ design. The precise
clique matchings (and partitions into precise clique matchings) with
$n=2^{t+1}$ and $q=2^t$ are constructed in \cite{41:Pot}.

 Mills in
 \cite{41:Mil} showed that for $n>3$, $n\neq 5$ an $H(n,q,4,3)$
design exists if and only if $nq$ is even and $q(n-1)(n-2)$ is
divisible by $3$. Ji in \cite{41:Ji} proved that an $H(5,q,4,3)$
exists if $q$ is even, $g\neq2$, and $q\not \equiv 10,26({\rm
mod48})$.

 Consider an
$H(n,q,w,t)$ design as a constant-weight code. The Hamming
distance\footnote{\,Here we consider elements of H-design as words
in alphabet $\{*,0,\dots,q-1\}$.} between two codewords of an
H-design is always greater than $w-t$.  The {\it code distance} of a
design is the  minimum Hamming distance between two codewords of
this design. The code distance of $H(n,q,w,t)$ design is at most
$2(w-t+1)$. An $H(n,q,w,t)$ design that forms a code with minimum
Hamming distance  $2(w-t+1)$ is called a generalized Steiner system
(see \cite{41:Et}). Note that an ordinary Steiner system
($H(n,1,w,t)$ design) is always a~code with Hamming distance
$2(w-t+1)$.

Etzion in \cite{41:Et} obtained some series of constructions of
generalized Steiner systems that are $H(n,q,3,2)$ or $H(n,2,4,3)$
designs. He proved that a generalized Steiner system being  an
$H(n,2,3,2)$ design exists if and only if $n\equiv 0\ \mbox {or}\ 1
(\mbox{mod}\,3)$, $n\geq 4$, $n\neq 6$.

Similarly we can consider an $A(n,q,w,t)$ design with the maximum
code distance. The code distance of an A-design is at most
$1+2(w-t)$ (but it can be equal to $1$). $A(n,2,n,n-1)$ designs with
Hamming distance $2$ were firstly constructed in \cite{41:Ham} for
every $n\geq 4$.
 Krotov
\cite{41:Kr} and Svanstr\"om \cite{41:Sv} proved (in other terms)
that $A(n,2,n,n-1)$ designs with  Hamming distance $3$ exist if and
only if $n=2^t$. It is straightforward that each $A(n,2,n,t)$ design
with Hamming distance $1+2(n-t)$ is a perfect ternary
constant-weight code.

\section{Constructions}

In this sections we consider some constructions of A- and H-designs
and partitions of the set of $m$-faces into A- and H-designs. Denote
by $Q^n_q(w)$ the set of $(n-w)$-faces of $Q^n_q$. Obviously
$|Q^n_q(w)|=q^w{n \choose w}$. It is easy to calculate that
cardinalities of $A(n,q,w,t)$ and $H(n,q,w,t)$ are equal to
$\alpha(n,q,w,t)=q^t\frac{n!(w-t)!}{w!(n-t)!}$. But the cardinality
of partition of the set $Q^n_q(w)$ into $H(n,q,w,t)$ designs is
equal to ${{n-t}\choose{n-w}}q^{w-t}=q^w {n\choose w}/
\alpha(n,q,w,t)$ and the cardinality of partition of the set
$Q^n_q(t)$ into $A(n,q,w,t)$ designs is equal to ${w\choose t}=q^t
{n\choose t}/ \alpha(n,q,w,t)$.

We propose the following constructions of H-designs.

{\bf Construction I}. Let  $S\subset Q_{q*}^n$ be an $H(n,q,w,t)$
design and let  $R\subset Q_{q'*}^w$ be an $H(w,q',w,t)$ design (MDS
code). Given $(a^1,\dots,*,\dots,a^i,\dots,*,\dots,a^w)\in S$ and
$(b_1,\dots,b_w)\in R$ arrange the  codeword\\
$((a^1,b_1),\dots,*,\dots,(a^i,b_i), \dots,*,\dots,(a^w,b_w))\in
Q_{qq'*}^n$. Let $T$ be the set of all these codewords.

\begin{proposition} \label{41:p:1}
 $T$ is an $H(n,qq',w,t)$ design.
\end{proposition}

{\it Proof.} Take $c^i\in Q_q$ and $d^i\in Q_{q'}$, and let
$((c^1,d^1), \dots,*,\dots,(c^i,d^i),\dots,*,$ $\dots,(c^t,d^t))$ be
arbitrary elements of $Q_{qq'*}^n$ with  weight $t$.  By the
definition of H-design there exists  a unique codeword
$(a^1,\dots,*,\dots,a^i,\dots,*,\dots,a^w)\in S$ such that the
$(n-w)$-face $(a^1,\dots,*,\dots,a^i,\dots,*,\dots,a^w)$ is
contained in the $(n-t)$-face $
(c^1,\dots,*,\dots,c^i,\dots,*,\dots,c^t)$. Convert the codeword\\
$(d^1,\dots,*,\dots,d^i,\dots,*,\dots,d^t)\in Q_{q'*}^n$ to the new
word with length $w$ removing a position $i$ if
$(a^1,\dots,*,\dots,a^i,\dots,*,\dots,a^w)$ has   $*$ in position
$i$. So, we form a codeword $\overline{d}\in Q_{q'*}^w$.
 By the definition of
H-design there exists a unique  codeword $(b_1,\dots,b_w)\in R$ such
that
 $(b_1,\dots,b_w)\subset \overline{d}$. Then the set $T$ is an $H(n,qq',w,t)$ design
 by  definition. $\bigtriangleup$

If we have partitions of the sets $Q^n_q(w)$ and $Q^w_{q'}(w)$ into
$H(n,q,w,t)$ and $H(w,q',w,t)$ designs respectively then  we obtain
a partition of the set $Q^n_{qq'}(w)$ into $H(n,qq',w,t)$ designs by
using Construction I for every pairs of $H(n,q,w,t)$ and
$H(w,q',w,t)$ designs from this partitions.

As mentioned above, $H(2k,k,2k-1,2k-2)$ designs exist for $k=2^t$,
$t\geq 1$. Since MDS codes with  distance $2$ ($H(m,q,m,m-1)$
designs) exist for all $q\geq 2$ and $m\geq 2$,  we get

\begin{corollary} \label{41:c:2}
For all $s,t\geq 1$ there exist
$H(2^{t+1},s2^t,2^{t+1}-1,2^{t+1}-2)$ designs.
\end{corollary}

Since partition of the sets $Q^{2^{t+1}}_{2^t}(2^{t+1}-1)$ into
$H(2^{t+1},2^t,2^{t+1}-1,2^{t+1}-2)$ designs exists \cite{41:Pot} it
is possible to construct  a partition of the set
$Q^{2^{t+1}}_{s2^t}(2^{t+1}-1)$ into
$H(2^{t+1},s2^t,2^{t+1}-1,2^{t+1}-2)$ designs for all $s,t\geq 1$.

 Note that  the MDS code $R$ in
Construction I can be chosen independently for every codeword from
$S$. The number of different $H(m,3,m,m-1)$ designs is $3\times
2^{m-1}$ (see \cite{41:KP}).
 A doubly
exponential lower bound of the number of MDS codes with distance $2$
($q\geq 4$) was established in \cite{41:KP}. Thus the number of
$H(2^{t+1},s2^t,2^{t+1}-1,2^{t+1}-2)$ designs is double exponential
with respect to the dimension $2^{t+1}$ as $s\geq 3$.

{\bf Construction II}. Let  $S\subset Q_{q*}^n$ be an $A(n,q,w,t)$
design.  For each pair of
$(a^1,\dots,*,\dots,a^i,\dots,*,\dots,a^t)\in S$ and
$(b_1,\dots,b_t)\in Q^t_{q'}$ we form the  codeword
$((a^1,b_1),\dots,*,\dots,(a^i,b_i),\dots,*,\dots,(a^t,b_t))\in
Q_{qq'*}^n$. Let $U$ be the set of all these codewords.

\begin{proposition} \label{41:p:3}
 $U$ is an $A(n,qq',w,t)$ design.
\end{proposition}
The proof  is similar to that of Proposition \ref{41:p:1}.

 As
mentioned above, each Steiner system $S(n-w,n-t,n)$ is equivalent to
an $A(n,1,w,t)$ design.

\begin{corollary} \label{41:c:3}
If there exists a Steiner system $S(n-w,n-t,n)$ then for each $q\geq
1$ there exists an $A(n,q,w,t)$ design.
\end{corollary}

It is easy to construct a partition of the set $Q^n_q(t)$ into
$A(n,q,w,t)$ designs from a partition of the layer of Boolean
$n$-dimensional cube into Steiner systems $S(n-w,n-t,n)$.

{\bf Construction III}. Let  $S\subset Q_{q*}^n$ be an
$A(n,q,n-1,n-2)$ design. Define $V=(S\times Q^n_q)\cup (Q^n_q\times
S)$.

\begin{proposition} \label{41:p:4}
 $V$ is an $A(2n,q,2n-1,2n-2)$ design.
\end{proposition}
{\it Proof.} Suppose that
$(c_1,\dots,c_{i-1},*,c_{i+1},\dots,c_{2n})$ is a word of weight
$2n-1$. If $i\leq n$ then there exists a unique codeword
$\overline{a}\in S$ such that\\
$(c_1,\dots,c_{i-1},*,c_{i+1},\dots,c_{n})\subset \overline{a}$. It
is clear that $(c_1,\dots,c_{i-1},*,c_{i+1},\dots,c_{2n})\subset
(\overline{a},\overline{d})$ where
$\overline{d}=(c_{n+1},\dots,c_{2n})$. The case $n<i\leq 2n$ is
similar.$\bigtriangleup$

\section{Multidimensional permanent}

  In \cite{41:Avg}  Avgustinovich developed a method of counting the
number of combinatorial configurations in terms of the
multidimensional permanent. Consider a biregular bipartite  graph
$G=(L,R,E)$ with parts $L$ and $R$. A set $C\subseteq L$ is called
$(L,R)${\it -perfect code}  if for each $v\in R$ there exist only
one vertex $u\in C$ such that $u$ is adjacent to $v$. The definition
of $(R,L)$-perfect code is obtained by changing parts $L$ and $R$.
It is easy to see that cardinalities of any $(L,R)${-perfect code}
and any $(R,L)$-perfect code of the same biregular bipartite  graph
are coincide.

Suppose that $\{ C_1,\dots,C_k\}$ is a partition of $L$ into
$(L,R)$-perfect codes. We define the {\it adjacency array}
$M(G,L)=(m_{i_1\dots i_k})$ by the following equation
 $m_{i_1\dots i_k}=|B^1_{i_1}\cap\dots \cap B^k_{i_k}|$ where
 $B^j_{i_j}$ is a neighborhood of the $i_j$th vertex of $C_j$.
If there exist a partition of $R$ consisted of $(R,L)$-perfect codes
then it is possible to define an adjacency array
$M(G,R)=(m_{i_1\dots i_k})$ by analogous way.

   A
$k$-element subset $I$ of $\{1,\dots,N\}^k$ is called a {\it
diagonal} if every pair of vectors  $\overline{i}, \overline{j}\in
I$ is distinct in each position that is $i_\sigma\neq j_\sigma$ for
all $\sigma\in\{1,\dots,k\}$. We define the $k$-{\it dimensional}
permanent of $M(G,L)$ as
$${\rm  per}_k M({ G,L} ) = \sum_{I \in D_N}\prod_{(i_1,\dots,i_k)\in I} m_{i_1\dots i_k},$$
where $D_N$ is the set of all diagonals. The following statement is
straightforward.

\begin{proposition} \label{41:p:15}
The number of $(R,L)$-perfect codes of $G$ is equal to ${\rm per}_k
M({ G,L} )$.
\end{proposition}

 Consider a $k$-partite
hypergraph $G_k$ containing $N$ vertices   in each part $C_i$,
$i=1,\dots,k$. Suppose that each $k$-edge of $G_k$ consists of $k$
vertices, with one vertex in each part of the hypergraph. A set of
disjoint $k$-edges that matches all vertices of the hypergraph is
called a {\it perfect $k$-matching}. Let each part $C_i$ of the
hypergraph be enumerated by $1,2,\dots,N$. We define the {\it
adjacency array} $M(G_k)=(m_{i_1\dots i_k})$ by the following rule:
$m_{i_1\dots i_k}=1$ if there exists a $k$-edge consisting of
vertices with numbers $i_1$ from the first part, $i_2$ from the
second part and so on and $m_{i_1\dots i_k}=0$  otherwise.

It is well known that the permanent of the adjacency matrix of a
bipartite graph is equal to the number of perfect matchings of the
graph.  The following statement is straightforward.

\begin{proposition} \label{41:p:5}
The number of perfect $k$-matchings of a hypergraph $G_k$ is equal
to ${\rm  per}_k M({ G_k} )$.
\end{proposition}

It is clear that any biregular bipartite  graph $G=(L,R,E)$ with
partition $\{ C_1,\dots,C_k\}$ of part $L$ into $(L,R)$-perfect
codes is equivalent to a $k$-partite regular hypergraph $G_k$ with
parts $C_1,\dots,C_k$. Here  $k$-edges of $G_k$ correspond to
vertices of the second part $R$ of $G$ and  perfect $k$-matchings of
$G_k$ one-to-one correspond to $(R,L)$-perfect codes of $G$.

 Given integers $w$, $t$ $(n\geq w> t\geq 1)$, define the bipartite
graph $G(n,q,w,t)$ with the parts $L=Q^n_q(w)$ and $R=Q^n_q(t)$. The
pair of vertices $\overline{c}\in Q^n_q(w)$ and $\overline{b}\in
Q^n_q(t)$ are connected by an edge in $G(n,q,w,t)$ if and only if
$\overline{c}\subset \overline{b}$. By  definition each $H(n,q,w,t)$
design is a subset of $Q^n_q(w)$ such that the neighborhoods of its
vertices do not intersect but cover
 $Q^n_q(t)$. We assume that
there exists a partition $H=\{H_1,\dots,H_k\}$, where
$k={{n-t}\choose{n-w}}q^{w-t}$, of  $Q^n_q(w)$ into $H(n,q,w,t)$
designs.

%Define the $k$-partite  hypergraph $GH$ with parts $H_1,\dots,H_k$.
%A collection $\{\overline{c}_1, \dots,\overline{c}_k\}$, where
%$\overline{c}_i\in Q^n_q(w)$, is a $k$-edge in $GH$ if there exists
%$\overline{b}\in Q^n_q(t)$, $\overline{b}=\bigcup_{i=1}^k
%\overline{c}_i$.

\begin{proposition} \label{41:p:6}
The number of different $A(n,q,w,t)$ designs is equal to\\ ${\rm
per}_k M({G(n,q,w,t),L} )$.
\end{proposition}{\it Proof.}
%Every $b\in  Q^n_q(t)$ contains $k$ different $(n-w)$-faces from
%$Q^n_q(w)$. By the definition of H-design all this $(n-w)$-faces lie
%in the different H-designs $H_1,\dots,H_k$.  Then the set $Q_q^n(t)$
%is in one-to-one correspondence to the set of $k$-edges in $GH$.
Any A-design $B\subset Q_q^n(t)$ perfectly covers all $(n-w)$-faces.
Then $B$ is a $(R,L)$-perfect code of $G(n,q,w,t)$. Using
Proposition \ref{41:p:15}, we obtain that the number of different
$A(n,q,w,t)$ designs is equal to ${\rm per}_k M({G(n,q,w,t),L} )$.
$\bigtriangleup$

% if there exists a partition $A=\{A_1,\dots,A_m\}$,
%where $m={{w}\choose{t}}$, of  $Q^n_q(t)$ into $A(n,q,w,t)$ designs,
%then we can define the $m$-part hypergraph $GA$ and state the
%following

 By  the definition each $A(n,q,w,t)$
design is a subset of $Q^n_q(t)$ such that  its faces do not
intersect and cover
 $Q^n_q(w)$. Let us to assume that
there exists a partition $A=\{A_1,\dots,A_m\}$, where $m={w\choose
t}$, of  $Q^n_q(t)$ into $A(n,q,w,t)$ designs.

%Define the $m$-partite hypergraph $GA$ with parts $A_1,\dots,A_m$. A
%collection $\{\overline{b}_1, \dots,\overline{b}_m\}$, where
%$\overline{b}_i\in Q^n_q(t)$, is an $m$-edge in $GA$ if there exists
%$\overline{c}\in Q^n_q(w)$, $\overline{c}=\bigcap_{i=1}^m
%\overline{b}_i$.
 Analogously to Proposition \ref{41:p:6}, we can
prove the following

\begin{proposition} \label{41:p:7}
The number of different $H(n,q,w,t)$ designs is equal to\\ ${\rm
per}_m M({G(n,q,w,t),R} )$.
\end{proposition}

The constructions of Section 2  provide examples of the parameters
such that there exists a partition into H-designs (A-designs). Thus
we can calculate the numbers of H-designs or A-designs with these
parameters  using the multidimensional permanent.

As mentioned above  $H(n,q,n,t)$ designs coincide with MDS codes
with distance $d=n-t+1$. For any integers $n,q,t$ there exists a
partition of $Q^n_q(t)$ into $A(n,q,n,t)$ designs where every
A-design consists of all parallel $t$-faces. Then the problem of
existence MDS codes with distance $d=n-t+1$ in $Q^n_q$ is equivalent
to the inequality ${\rm per}_m M(G(n,q,w,t),R )>0$. Ball
\cite{41:Ball} proved that linear MDS code over prime field
$\mathbb{F}_q$ has length at most $q+1$ (except for the trivial
cases $d=2,n$). But the question of existence of nonlinear MDS codes
of larger lengths is open.

\section{Acknowledgements} The author thanks D. S. Krotov for his
interest in this work.

\bibliographystyle{model1a-num-names}

\end{document}